# REJOINDER

By Marc Hallin, Davy Paindaveine and Miroslav Šiman

*Université Libre de Bruxelles*

We heartily thank all discussants for their thoughtful and friendly comments on our paper. As the three discussions address distinct issues, with very little overlap, this rejoinder is also organized into three distinct parts.

## 1. Discussion by Serfling and Zuo.

### 1.1. *Point-valued versus hyperplane-valued quantiles.*
Univariate quantiles are point-valued. Points on the real line $\mathbb{R}^1$ are hyperplanes of dimension 0, and hence, in a $k$-dimensional context, defining quantiles as points or as $(k-1)$-dimensional hyperplanes may be equally legitimate and, depending on the properties to be emphasized or the inference problem at hand, equally sensible and useful.

If $L$-estimation is to be privileged, point-valued quantiles, which provide easily tractable integrands for $L$-functionals, may look more attractive. Many basic properties of quantiles, however, are better accounted for by hyperplane-valued generalizations. In particular, the role of univariate quantiles as critical values partitioning, with respect to outlyingness, the sample or the observation space into two regions with given frequencies or probability contents cannot be assumed, in $\mathbb{R}^k$, by any point-valued concept. Since our objective was to establish a meaningful and fully operational bridge between the quantile and depth universes, hyperplane-valued concepts quite naturally came into the picture; this does not imply, however, that point-valued quantiles should be abandoned in favor of hyperplane-valued ones.

### 1.2. *Robustness issues: A useful adjunct to the halfspace DOQR paradigm.*
The discussants point out that our quantile hyperplane construction, since it leads to halfspace depth, also inherits all the properties of the halfspace depth DOQR (Depth, Outlyingness, Quantile, Rank) paradigm. Among them, they are stressing [point (H2) of their discussion] a poor performance in outlier detection.

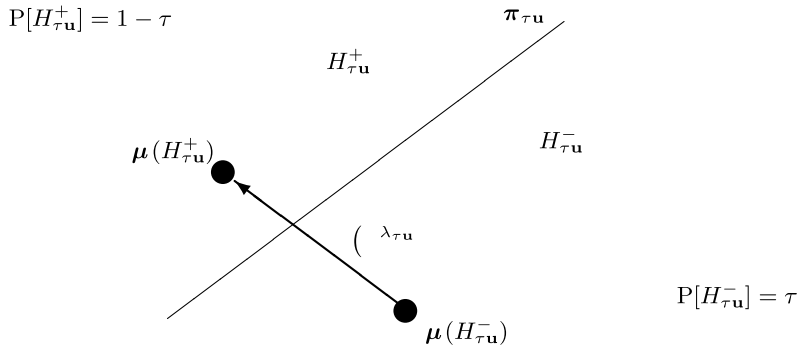

Fig. 1.  *An illustration of the mutual relationship between $\tau$, $\mathbf{u}$ and $\lambda_{\tau\mathbf{u}}$ (for appropriate $\tau$ values); $\boldsymbol{\mu}(H_{\tau\mathbf{u}}^+)$ and $\boldsymbol{\mu}(H_{\tau\mathbf{u}}^-)$ stand for the probability mass centers of $H_{\tau\mathbf{u}}^+$ and $H_{\tau\mathbf{u}}^-$, respectively.*

Since no robustness requirement has been made at any stage in our construction of directional quantiles, such poor performance is not surprising. Our quantile hyperplanes, however, do not come *alone*: thanks to their $L_1$ nature, and their relation to linear programming, each $\pi_{\tau\mathbf{u}}^{(n)}$ indeed is accompanied by a Lagrange multiplier value $\lambda_{\tau\mathbf{u}}^{(n)}$, which is extremely helpful for outlier identification, as we now explain. While, for any $\tau \in (0,1)$, the halfspace depth contour of order $\tau$ is not affected by moving its interior points, the Lagrange multiplier process $\boldsymbol{\tau} := \tau\mathbf{u} \mapsto \lambda_{\boldsymbol{\tau}}^{(n)}$ is extremely sensitive to such changes, hence clearly contains useful information which is not directly available from the contour itself.

More specifically, the mutual relationship between $\tau$, $\mathbf{u}$ and $\lambda_{\tau\mathbf{u}}^{(n)}$—see, in the population version, equation (3.6) or Figure 1—makes clear that large values of $\lambda_{\tau\mathbf{u}}^{(n)}$ indicate that the probability mass centers associated with observations lying above and below $\pi_{\tau\mathbf{u}}^{(n)}$ are far apart from each other, which may reveal the presence of outliers (or some clustering structure): see Figure 2 for an empirical illustration. Incidentally, this multiplier was shown in [4] to be closely related, in a portfolio optimization context, to the *shortfall risk measure* or *Conditional VaR* (CVaR), which is widely used in finance for detection of extreme observations. The quantile hyperplane construction of depth contours thus complements the contours themselves with a most useful outlier-detection tool: the $\lambda_{\tau\mathbf{u}}^{(n)}$'s.

1.3.  *Miscellaneous issues and questions.*  We finally comment on various other points raised throughout the discussion.

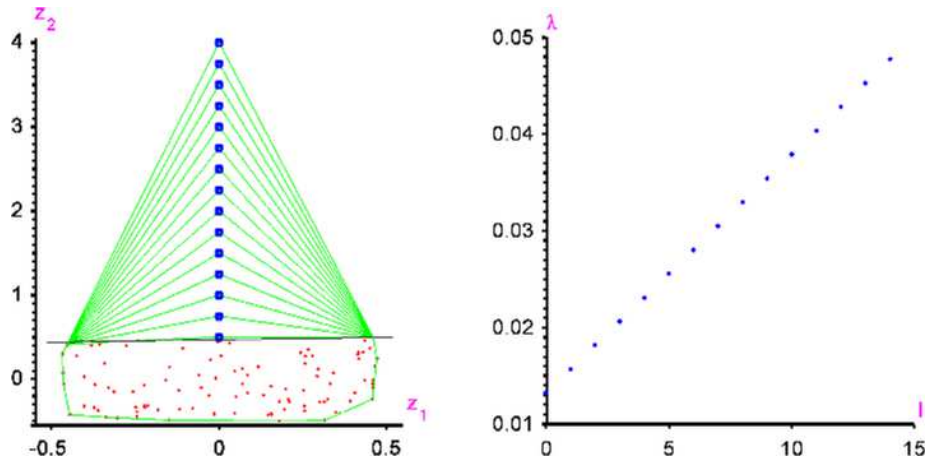

Fig. 2. *The left plot shows $n - 1 = 98$ (red) points drawn from $U([-0.5, 0.5]^2)'$, the centered bivariate uniform distribution over the unit square, and one (blue) outlier of the form $(0, 0.5 + \ell/4)'$, $\ell = 0, \ldots, 14$, along with, for all $\ell$, (i) the halfspace depth contour of order $1/n$ (in green), that is, the convex hull of the sample, and (ii) the $(\tau \mathbf{u}_0)$-quantile hyperplane (in black) obtained for $\tau = 2.5/n$ and $\mathbf{u}_0 = (0, -1)'$. The right plot reports values of the corresponding Lagrange multiplier $\lambda_{\tau \mathbf{u}_0}^{(n)}$ for each of the 15 samples considered. The halfspace depth region of order $1/n$ clearly breaks down as expected, whereas the $(\tau \mathbf{u}_0)$-quantile hyperplane is not affected by the position of the blue point. Still, the outlyingness of that blue point is revealed in the corresponding Lagrange multiplier value.*

*Smoothness of contours.* The discussants feel [requirement (C8) in their discussion] that polygonal, hence nonsmooth contours, are a weakness. Still, the more complicated a contour, the more difficult its description, handling, plotting and exact quantitative analysis, especially in high-dimensional spaces. In particular, determining whether a given point is lying inside or outside a given contour, or computing a contour's (estimated) probability contents, which may be of primary practical importance, may be much more difficult for smooth contours than for polyhedral ones: smoothness then actually is a drawback.

*Connections with univariate quantiles and quantile analysis.* We are sorry for the misunderstanding caused by our explanation. The vector representation of $\boldsymbol{\tau}$ was used to stress that the sign, in the univariate case, also matters: in the population case, our lower quantile halfspaces reduce to $(-\infty, F^{-1}(\tau))$ and $(F^{-1}(1-\tau), +\infty)$ for $\boldsymbol{\tau} = \tau.1$ and $\boldsymbol{\tau} = \tau.(-1)$, respectively, each of these halfspaces carrying a probability mass equal to $\tau$, as announced in (3.2a). Finally, we point out that our $\boldsymbol{\tau}$-quantiles are well defined for any $\tau \in (0, 1)$, even in the sample case, because $N = 0$ is also a possibility there. So (3.9) does not restrict the range of $\tau$ in any way; nor does it in the standard single-output quantile regression framework.

*Moment and regularity assumptions.* We agree that the assumptions are relatively strong (not quite as strong, though, as those required for many other concepts). What we would like to point out is our remark at the end of Section 3.1 (more elaborated on in [4]) that the results in the empirical case [such as (3.11)] were derived for any dataset in general position; as for (3.8), (3.9) and (5.2), they do not require even that. So we do not need finite moments, independence, continuity or unimodality in the analysis of empirical data. In other words, our methodology *always* keeps the frequency of directional outliers ($N/n$) under control and is thus as nonparametric as possible in that respect. This control even extends to the local case considered in [1]. This allows us to consider arbitrary regressors in models with arbitrary error terms. And if the frequency $N/n$ asymptotically differs from the outlying probability $P[\mathbf{Z} \in H_\tau^-]$, then we feel that the former should be favored, which is what most technical norms already do (of course, asymptotic theory would still be required for inference).

*The assumption of "general position."* (i) Although we were able (in [1]) to relax this assumption to allow for multiple identical observations (which makes resampling methods technically possible), we would indeed very much like to know whether that assumption can be totally removed, as this would allow, for example, direct inclusion of discrete or dummy variables in the regression.[1] Presently, one has to investigate the model for each dummy value separately (which seems nevertheless quite reasonable) and possibly also perturb highly discrete data with random noise of reasonably small magnitude. (ii) Indeed, the fraction $\ell/n$ must be lower than the maximum halfspace depth to keep the assumption of Theorem 4.2 satisfied. We do not see this as a severe restriction because we still establish the contour equivalence there for all nondegenerate halfspace depth contours. This does not imply that a stronger statement does not hold.

*Compactness of the $R(\tau)$ regions.* "Empty" here really means "empty" and not "not containing any sample point," which might be the cause of confusion. With this clarified, the polyhedral regions mentioned by the discussants do not contradict our claims.

*Usefulness of our quantile hyperplanes.* Finally, we stress that, beyond their role in outlier identification (see Section 1.2 above), our quantile hyperplanes might also be used to define a location-free concept of angular distance. More precisely, one could modify the concept of interdirections from [6] by replacing the separating hyperplanes passing through the origin there by our 0.5-quantile hyperplanes. This would probably provide a meaningful measure of angular distance between two observations ("detrended," in the regression case), at least for (conditional) halfspace symmetric distributions and with no need of data centering.

---

[1] Bad configuration of data points can indeed make exact computation of the contours by our current codes impossible; further improvement in that respect would be welcome.

## 2. Discussion by Wei.

2.1. *Generality, flexibility and nonparametric extensions of the method.* The goal of our Section 6 was to show that the multivariate quantiles proposed in Section 2 are easily extended into multiple-output *regression* quantiles. This is achieved, quite simply, by restricting to those quantile hyperplanes associated with directions $\mathbf{u} \in \mathbb{R}^k$ in the subspace of the response. As a direct consequence, our multiple-output regression quantiles still are hyperplanes, hence linear objects. Although no assumption on the actual regression structure (i.e., the relation between $\mathbf{x}$ and the distribution of $\mathbf{Y}$ conditional on $\mathbf{X} = \mathbf{x}$) is required in order to derive consistency and other asymptotic results (see our Section 3.3), such linear objects cannot be expected to provide a satisfactory fit and a meaningful (local or conditional) interpretation unless the underlying regression structure itself is linear.

We thus fully agree with the discussant that a nonparametric extension of our method, able to accomodate nonlinear structures, would be welcome. We would tend to avoid, however, the piecewise polynomial regression approach suggested in this discussion, which involves too many unknown parameters to be set a priori. For one random regressor ($p = 2$), this includes the number of knots and their positions, the degrees of the polynomials to be considered on each interval, and various smoothness constraints at the knots and both end intervals. For more than one random regressor ($p \geq 3$), such choices become untractable. A further disadvantage of the piecewise polynomial approach lies in the fact that each additional term in the polynomials increases by one unit the total number of regressors, which slows down computation, while producing less stable final results.

In [1], we therefore rather suggest either a local constant on a local linear approach, on the model of the local linear single-output regression quantiles proposed by [3] and [7]. Such an approach fits nicely into a conditional halfspace depth context, and is moreover able to capture complicated forms of heteroskedasticity. Of course, the choice of an adequate (possibly matrix-valued) bandwidth remains a delicate point, but this is probably unavoidable in a setting involving ($p - 1$) random regressors. When applied to the body girth measurement dataset, that local linear method seems to perform very well, and takes care of the drawbacks of the *global linear* approach the discussant is rightly pointing out. We refer to [1] for details.

2.2. *Diagnostics.* If the statistic $\Delta(\mathbf{x})$ proposed by the discussant were computed on the basis of the available data points, it would always be roughly equal to $\tau$ [this follows from our empirical subgradient conditions, which are satisfied in the absence of any assumptions on the underlying distribution of $(\mathbf{X}, \mathbf{Y})$]. Therefore, we rather see $\Delta(\mathbf{x})$, or its integral over the regressor space, as a measure of the *nonlinearity* of the regression, hence as an interesting test statistic for the null hypothesis of linearity. Equivalently, such statistic measures how the global quantile behavior differs from the local one, which might be of interest since our global regression regions have no clear link to the conditional distributions associated with given regressor values.

Although we suggested to use $\lambda_\tau$ or $\|\mathbf{c}_\tau\|$ for testing symmetry as an illustration of possible applications in statistical inference, this by no means was an exhaustive list of examples. For instance, the process $\boldsymbol{\tau} \mapsto \lambda_{\boldsymbol{\tau}}^{(n)}$ can be used for elliptical goodness-of-fit tests as well, since its population counterpart $\boldsymbol{\tau} \mapsto \lambda_{\boldsymbol{\tau}}$ is known exactly for elliptical distributions; see [4]. Besides, if one works with the constrained optimization version of Definition 6.1 (i.e., the one obtained by extending Definition 2.2 rather than Definition 2.1 to the regression setup), one could consider the subvectors of $\mathbf{c}_{\boldsymbol{\tau}}$ related to those regressors whose joint importance is to be tested.

## 3. Discussion by Kong and Mizera.

3.1. *"Quantile or not quantile?"* Since 1978 and Koenker and Bassett's epochal contribution, statisticians have been calling *quantile hyperplane* any hyperplane dividing a data cloud into two parts with given proportions of observations below and above and minimizing some $L_1$ criterion (for simplicity, we restrict ourselves to the empirical case). This is exactly what "our" $\pi_{\mathrm{HP\check{S}};\tau\mathbf{u}}^{(n)} := \pi_{\tau\mathbf{u}}^{(n)}$ hyperplanes are doing; this is also what the $\pi_{\mathrm{KM};\tau\mathbf{u}}^{(n)}$ hyperplanes from Kong and Mizera [2] are doing, the only (yet fundamental) difference being the choice of the $L_1$ criterion (before or after projection) to be minimized or, equivalently, the class of hyperplanes (unrestricted or orthogonal to $\mathbf{u}$) over which minimization is performed. Therefore, the quantile terminology here is not only nonusurped, but also highly adequate.

While Kong and Mizera refrain from emphasizing those hyperplanes, refrain from providing them with a specific notation, and refrain from calling them quantiles, those $\pi_{\mathrm{KM};\tau\mathbf{u}}^{(n)}$'s are pervasively present throughout their text and play a major role in the definition and interpretation of the envelopes that eventually coincide with halfspace depth contours.

Actually, quantile ideas are so much at the core of the problem under discussion that this avoidance of the quantile terminology is verging the oulipian exercise of liponomy.[2] *Quantile or not quantile*, thus, *that is* not *the question*. Not here. The question is about two different quantile-based ways of characterizing contours that happen to be the halfspace depth contours in the location case. Our paper opens with a marked homage to [2], and we keep on thinking that the relation established there between univariate directional quantiles (equivalently, the infinite collection of $\pi^{(n)}_{\mathrm{KM};\tau\mathbf{u}}$ hyperplanes) and depth contours is a remarkable and most important result. Nevertheless, with probability one, and avoiding integer values of $n\tau$, the exact halfspace depth contours are obtained via the *finite* collection of $\pi^{(n)}_{\mathrm{HPS};\tau\mathbf{u}}$ hyperplanes, which is a finite subset[3] of the infinite collection of $\pi^{(n)}_{\mathrm{KM};\tau\mathbf{u}}$'s. We therefore think that our construction is more direct, more parcimonious and—see Section 3.2 below—equally meaningful.

3.2. *Talking to health professionals.* We definitely share the discussants' concern about being able to explain statistical concepts and analyses to people from other fields. This concern is fundamental to our profession (which does not imply that *The Annals of Statistics* should be written for health professionals).

Going back to the BMI example developed in the discussion, we do not see the slightest pedagogy problem there. For simplicity, let us start with $\tau = 0.5$, that is, with medians. Once the axes have been rotated as in Figure 1 of the discussion, the $\pi^{(n)}_{\mathrm{KM};0.5\mathbf{u}}$ ($\mathbf{u} = (0,1)'$) hyperplane is the horizontal hyperplane through the univariate median of $\log(\mathrm{BMI})$—in other words, the *horizontal* hyperplane that fits the data cloud best from the point of view of absolute $\log(\mathrm{BMI})$ deviations. Similarly, the $\pi^{(n)}_{\mathrm{HPS};0.5\mathbf{u}}$ hyperplane is the (not necessarily *horizontal*) hyperplane that fits the data cloud best from the point of view of absolute $\log(\mathrm{BMI})$ deviations. Where is the intuition problem? Practitioners are perfectly grasping the meaning of LAD and quantile regression, which by now are part of daily practice. Switching to weighted positive and negative absolute deviations, the $\tau = 0.75$ case is hardly more difficult. Now, this is essentially a rephrasing of the discussants' allegedly unintuitive conditional quantile description.

As for the role of the direction $\mathbf{u}$, our paper provides a meaningful interpretation—illustrated in Figure 1. Besides, the link between hyperplanes and directions is less important once the full contours are considered.

---

[2]A stylistic device, which consists in writing without the use of certain words. Liponomy was popularized by the literary group of writers and mathematicians known as Oulipo (Ouvroir de Littérature Potentielle), of which Raymond Queneau, Italo Calvino and Georges Perec were influential members. Georges Perec himself is most famous for *La Disparition* (1969), a full 319-page novel avoiding the letter "e."

[3]After adequate $\mathbf{u}$-relabeling (for any $\mathbf{u}$, there exists a $\mathbf{v}$ such that $\pi^{(n)}_{\mathrm{HPS};\tau\mathbf{u}} = \pi^{(n)}_{\mathrm{KM};\tau\mathbf{v}}$).

3.3. *Approximation or exact computation?*  As already mentioned, Kong and Mizera [2] basically address the same problem as we do. While their approach at best provides excellent approximate solutions, ours is leading to exact ones.[4] Referring to the toy example displayed in Figure 4 of our paper, the polygonal contours we obtain from the full collections of (eight and 18, respectively) quantile hyperplanes for the two $\tau$ values considered there both consist of eight facets; Kong and Mizera's construction, if based on 21 or 201 equispaced directions, as proposed in Section 4 of their discussion, yields hopefully excellent approximations of those eight-faceted polygons by 21- or 201-faceted ones. While we agree that drawing a curve, as a rule, "consists in drawing one hundred or one thousand small segments," drawing straight lines, in which case one segment is sufficient, is a notable exception.

Now, in the bivariate location case, we fully admit that the Kong and Mizera strategy, based on, say, $K = 2001$ directions, can provide such a perfect approximation that eye-inspection hardly detects any deviation from the exact contour. But it also produces a contour with 2001 vertices and facets. Even if the exact contour has as many as 100 facets, computing its volume or checking whether a given point is inside or not[5] can be expected to be about 20 times faster than from a polygon with 2001 facets. Moreover, it is by no means clear that such an approximation always would be computationally faster than the exact solution. The point is that the $K$ univariate quantiles for the approximate solution must be computed from scratch whereas once the first initial solution in the parametric programming approach is obtained, the recursive computation of the subsequent ones is much faster.

Still, our approach will be defeated by the computational burden for large $k$ and/or $n$, while a sampling-the-directions[6] strategy—based on the hyperplanes $\pi_{\mathrm{KM};\tau\mathbf{u}}^{(n)}$ associated with a fixed number, $K = 20{,}001$, say, of directions $\mathbf{u}$—remains implementable, basically, for any $n$. An alternative approximate construction, however, can also be based on the $K$ $\pi_{\mathrm{HPŠ};\tau\mathbf{u}}^{(n)}$ hyperplanes associated with the same sample of $\mathbf{u}$'s, which leads to another approximate solution. The latter is based on $K$ hyperplanes from the finite[7] collection

---

[4]Of course, rounding off and other numerical limitations are present in all implementations: when we say "exact solution," we mean that our method is yielding, up to numerical approximations, a polyhedral solution with $K_0$ facets which coincides with a Tukey contour with $K_0$ facets.

[5]A problem of primary importance if depth contours are to be interpreted as *quantile contours*, that is, reindexed by an estimation of their probability contents.

[6]"Sampling" here covers the systematic sampling method decribed in [2] but also possibly more sophisticated randomized or data-driven choices of the $K$ directions to be considered.

[7]Still avoiding the troublesome $\tau = i/n$ ($i \in \mathbb{N}$) values; moreover, unless $K$ is very small compared with $n$, the $\pi_{\mathrm{HPŠ};\tau\mathbf{u}}^{(n)}$ hyperplanes associated with neighboring sampled $\mathbf{u}$'s often

of $\pi^{(n)}_{\mathrm{HPŠ};\tau\mathbf{u}}$'s characterizing the actual halfspace depth contour, whereas the probability[8] that any of the $K$ hyperplanes $\pi^{(n)}_{\mathrm{KM};\tau\mathbf{u}}$ involved in the description of its Kong and Mizera counterpart contains any of the actual halfspace depth facets is strictly zero.

All comments above hold even in the bivariate location case, which is the most favorable one to Kong and Mizera [2]. As the dimension increases, the Kong and Mizera approximation clearly deteriorates (provided that the number of directions $K$ considered remains fixed), so that the dominance of our approach over theirs is even stronger. As for the multiple-output regression setting, we totally agree, as mentioned in our reply to Wei's discussion, that something should be done in order to capture possibly nonlinear dependencies in a flexible way. Such an extension of our approach was announced in Section 8 and is actually fully developed in [1]. We gratefully welcome the interesting references provided by the discussants, and their suggestions for overcoming the quantile crossing phenomenon.

3.4. *Implementation.* We definitely do not claim any expertise in algorithmics and computational geometry and we humbly apologize if our bibliography missed some major references. Although there are some papers addressing the notoriously hard problem of computing the halfspace *depth of a point* in any dimension, as the references proposed by the discussants also suggest, we are not aware of any exact *implementable* algorithm for computing halfspace depth *contours* beyond dimension two, which is quite a different problem. When stating this, we actually reproduce the information from recognized experts we have been consulting. In fact, as the discussants themselves admit, practical implementability is quite a serious issue even in the bivariate location case.

As far as numerical accuracy and stability of the implementation are concerned, these problems are not as serious as they might seem at first sight. Our solution (see [5]) requires only $k \times k$ matrices to be inverted, does not necessarily accumulate numerical errors in the course of computation, and although it requires an initial solution to start, it is only the identification of exactly fitted observations that really matters. When this identification seems unclear from the first interior-point solution obtained, any other solution corresponding to a different direction can be used instead. And although current implementation cannot handle degeneracy well, it could be used to detect it and to warn the user that the data should

---

coincide, allowing for much shorter contour descriptions. This is not the case for the $K$ $\pi^{(n)}_{\mathrm{KM};\tau\mathbf{u}}$'s, which are all distinct.

[8]Under absolute continuity.

be perturbed with a small random noise. This takes care of most problems of degeneracy or inaccurate initial solutions, even in the pathological linear case considered by the discussants, where a noise of amplitude as low as 0.00001 already does the job provided that the number of observations on the line segment is not too high. The troublesome values of $\tau$ leading to integer values of $n\tau$ must always be avoided, but this, in view of our Theorem 4.2, can be done without any loss of generality in the location case. The pictures in our paper testify that our implementation can handle at least ten thousand bivariate observations in a regression setting.

E.C.A.R.E.S.
UNIVERSITÉ LIBRE DE BRUXELLES
50, AVENUE F.D. ROOSEVELT, CP114
B-1050 BRUXELLES
BELGIUM
E-MAIL: mhallin@ulb.ac.be
        dpaindav@ulb.ac.be
        Miroslav.Siman@ulb.ac.be
URL: http://homepages.ulb.ac.be/~dpaindav